\documentclass[11pt]{amsart}
\usepackage{amssymb}
\usepackage{amscd}
\input amssym.def
\input amssym  
\newtheorem{theorem}{Theorem}[section]

\newtheorem{lemma}{Lemma}[section]

\newtheorem{remark}{Remark}[section]
\newtheorem{definition}{Definition}[section]
\newtheorem{proposition}{Proposition}[section]

\newcommand{\cal}{\mathcal}

\usepackage{amssymb}

\begin{document}
 
\title[Vertex-algebraic structure of principal subspaces]
{Vertex-algebraic structure of the principal subspaces of certain
$A_1^{(1)}$-modules, I: level one case}
%\footnotemark[\value{footnote}]

\author{C. Calinescu, J. Lepowsky and A. Milas}
%\footnotemark[\value{footnote}] 

\thanks{C.C. gratefully acknowledges partial
support {}from the Center for Discrete Mathematics and Theoretical
Computer Science (DIMACS), Rutgers University. J.L. gratefully
acknowledges partial support {}from NSF grant DMS--0401302.}

\begin{abstract}
This is the first in a series of papers in which we study
vertex-algebraic structure of Feigin-Stoyanovsky's principal subspaces
associated to standard modules for both untwisted and twisted affine
Lie algebras.  A key idea is to prove suitable presentations of
principal subspaces, without using bases or even ``small'' spanning
sets of these spaces.  In this paper we prove presentations of the
principal subspaces of the basic $A_1^{(1)}$-modules.  These
convenient presentations were previously used in work of
Capparelli-Lepowsky-Milas for the purpose of obtaining the classical
Rogers-Ramanujan recursion for the graded dimensions of the principal
subspaces.
\end{abstract}

\maketitle

\renewcommand{\theequation}{\thesection.\arabic{equation}}
\renewcommand{\thetheorem}{\thesection.\arabic{theorem}}
\setcounter{equation}{0}
\setcounter{theorem}{0}
\setcounter{section}{0}
\renewcommand{\theequation}{\thesection.\arabic{equation}}
\setcounter{equation}{0}

\section{Introduction}

One of the central problems in combinatorial representation theory is
to find (combinatorial) bases of standard modules for affine Kac-Moody
Lie algebras. Finding a basis of a standard module is closely related
to finding a basis of the vacuum subspace for an appropriate
Heisenberg algebra by means of vertex operators and ``$Z$-algebras''
(see in particular \cite{LW1}--\cite{LW4}, \cite{LP1}--\cite{LP2},
\cite{MP1}--\cite{MP2}).

In \cite{FS1}--\cite{FS2}, Feigin and Stoyanovsky associated to every
standard $A_n^{(1)}$-module $L(\Lambda)$ another distinguished
subspace, the ``principal subspace''
$$W(\Lambda)=U(\bar{\goth{n}}) \cdot v_{\Lambda},$$
where $\goth{n}$ is the nilradical of a Borel subalgebra of
$\goth{sl}(n+1)$,
$$\bar{\goth{n}}=\goth{n} \otimes \mathbb{C}[t,t^{-1}],$$ and
$\Lambda$ and $v_{\Lambda}$ are the highest weight and a highest
weight vector of $L(\Lambda)$.  Of course, this definition extends to
an arbitrary highest weight module for an affine Lie algebra.
Compared to the vacuum subspaces of $L(\Lambda)$ for Heisenberg Lie
algebras, the principal subspace $W(\Lambda)$ is a somewhat simpler
object, yet interestingly enough it still carries a wealth of
information about the standard module (see \cite{FS1}--\cite{FS3}).
{}From the definition, $W(\Lambda)$ can be identified with the quotient
$$U(\bar{\goth{n}})/I_{\Lambda},$$
where $I_{\Lambda}$ is the kernel of the natural map
$$f_{\Lambda} : U(\bar{\goth{n}}) \longrightarrow W(\Lambda), \ \ a
\mapsto a \cdot v_{\Lambda}.$$
It is an important general problem to find a presentation of
$W(\Lambda)$ (i.e., to precisely describe $I_{\Lambda}$), to compute
the corresponding graded dimension of $W(\Lambda)$, and especially,
{}from our point of view, to reveal new vertex-algebraic structure
underlying such issues.

A central aim of this paper is to supply a new proof of the natural
presentation of the principal subspaces $W(\Lambda_i)$ of the basic
$A_1^{(1)}$-modules $L(\Lambda_i)$, $i=0,1$, established in
\cite{FS1}, \cite{FS2}.  There are a number of proofs, all of which
have essentially used detailed structure, such as bases, of the
standard modules.  One wants to know how to prove this result using
less such structure.  For a recent proof in this direction, see
\cite{Cal}.  The presentation result was used in \cite{CLM1} as a step
in proving the classical Rogers-Ramanujan recursion for the graded
dimensions of the spaces $W(\Lambda_i)$.

The proof of the presentation that we provide in this paper does not
use any bases of the standard modules or of the principal subspaces,
or even ``small'' spanning sets that turn out to be bases (that is,
spanning sets whose monomials satisfy the difference-two condition on
their indices), and so gives new insight into the program of obtaining
Rogers-Ramanujan-type recursions initiated in \cite{CLM1} (for level
one modules) and \cite{CLM2} (for higher level modules).  In those
papers, exact sequences among principal subspaces were built, yielding
$q$-difference equations and in turn formulas for the graded
dimensions of the principal subspaces.  Vertex (operator) algebras and
intertwining operators among standard modules were used, together with
certain translations in the affine Weyl group of $\goth{sl}(2)$.  It
turns out that we have been able to use these same ingredients in a
new way, combined with the viewpoint of generalized Verma modules, to
prove the presentation result by induction on conformal weight.

Our main results are described in Theorems \ref{th1} and \ref{th2},
and in Theorem \ref{th3} we offer a reformulation of part of Theorem
\ref{th2} using the notion of ideal of a vertex algebra.

The ideas and arguments presented here certainly generalize, as we
will show in a series of forthcoming papers
\cite{CalLM1}--\cite{CalLM3}; the paper \cite{CalLM1} covers the case
of the principal subspaces of the higher level standard
$\widehat{\goth{sl}(2)}$-modules, for which the Rogers-Selberg
recursions were obtained in \cite{CLM2}.

\section{Formulations of the main result }

We use the setting of \cite{CLM1} (cf. \cite{FLM}, \cite{LL}).  We
shall work with the complex Lie algebra
$$
\goth{g}=\goth{sl}(2)= \mathbb{C}x_{-\alpha} \oplus \mathbb{C} h
\oplus \mathbb{C} x_{\alpha},
$$
with the brackets
$$
[h, x_{\alpha}]= 2 x_{\alpha}, \; \; [h, x_{-\alpha}]= -2
x_{-\alpha}, \; \; [x_{\alpha}, x_{-\alpha}]= h.
$$
Take the Cartan subalgebra $\goth{h}=\mathbb{C} h$.  The standard
symmetric invariant bilinear form $\langle a, b \rangle= \mbox{tr} \;
(ab)$ for $a,b \in \goth{g}$ allows us to identify $\goth{h}$ with
its dual $\goth{h}^{*}$. Take $\alpha$ to be the (positive) simple root
corresponding to the root vector $x_{\alpha}$. Then under our
identifications,
$$h= \alpha$$ and
$$ \langle \alpha, \alpha \rangle =2.$$

The corresponding untwisted affine Lie algebra is
\begin{equation}
\widehat{\goth{g}}= \goth{g} \otimes \mathbb{C}[t, t^{-1}] \oplus
\mathbb{C} {\bf k},
\end{equation}
with the bracket relations
\begin{equation} \label{brackets}
[a \otimes t^m, b \otimes t^n]=[a,b] \otimes t^{m+n} + m \langle a, b
\rangle \delta_{m+n,0}{\bf k},
\end{equation}
for $a, b \in \goth{g}$, $m,n \in \mathbb{Z}$, together with the
condition that ${\bf k}$ is a nonzero central element of
$\widehat{\goth{g}}$.  Let
$$
\tilde{\goth{g}}=\widehat{\goth{g}} \oplus \mathbb{C}d,
$$
where $d$ acts as follows:
$$
[d, a \otimes t^n]=na\otimes t^n \; \; \mbox{for} \; \; a \in
\goth{g}, \; n \in \mathbb{Z},
$$
$$
[d, {\bf k}]=0.
$$

Setting
$$
{\goth{n}}=\mathbb{C}x_{\alpha}, 
$$
consider the following subalgebras of $\widehat{\goth{g}}$:
$$
\goth{g}= \goth{g} \otimes t^0,
$$
$$
\bar{\goth{n}}={\goth{n}} \otimes \mathbb{C}[t, t^{-1}],
$$
$$
\bar{\goth{n}}_{-}= {\goth{n}} \otimes t^{-1}\mathbb{C}[t^{-1}]
$$
and
$$
\bar{\goth{n}}_{\leq -2} ={\goth{n}} \otimes t^{-2}
\mathbb{C}[t^{-1}].
$$

The Lie algebra $\widehat{\goth{g}}$ has the triangular decomposition
\begin{equation}
\widehat{\goth{g}}=( \mathbb{C} x_{-\alpha} \oplus \goth{g} \otimes
t^{-1} \mathbb{C}[t^{-1}]) \oplus (\goth{h} \oplus \mathbb{C}{\bf k})
\oplus (\mathbb{C}x_{\alpha} \oplus \goth{g} \otimes t\mathbb{C}[t]).
\end{equation}
It also has the triangular decomposition
\begin{equation}
\widehat{\goth{g}}= \widehat{\goth{g}}_{<0} \oplus
\widehat{\goth{g}}_{\geq 0},
\end{equation}
where
$$
\widehat{\goth{g}}_{<0}= \goth{g} \otimes t^{-1} \mathbb{C} [t^{-1}]
$$
and
$$
\widehat{\goth{g}}_{\geq 0}= \goth{g} \otimes \mathbb{C}[t] \oplus
\mathbb{C}{\bf k}.
$$

We consider the level $1$ standard $\widehat{\goth{g}}$-modules
$L(\Lambda_0)$ and $L(\Lambda_1)$ (cf. \cite{K}), where $\Lambda_0,
\Lambda_1 \in (\goth{h} \oplus \mathbb{C}{\bf k})^*$ are the
fundamental weights of $\widehat{\goth{g}}$ ($\langle \Lambda_i, {\bf
k} \rangle =1$, $\langle \Lambda_i, h \rangle=\delta_{i,1}$ for
$i=0,1$). Denote by $v_{\Lambda_i}$ the highest weight vectors of
$L(\Lambda_i)$ used in \cite{CLM1}, namely,
\begin{equation} \label{v}
v_{\Lambda_0}=1 \; \; \mbox{and} \; \; v_{\Lambda_1} =e^{\alpha/2}
\cdot v_{\Lambda_0}.
\end{equation}
(See Section 2 of \cite{CLM1}.)

The principal subspaces $W(\Lambda_i)$ of $L(\Lambda_i)$ are defined by
\begin{equation} \label{ps}
W(\Lambda_i)= U(\bar{\goth{n}}) \cdot v_{\Lambda_i}
\end{equation}
for $i=0,1$ in [FS1]. 
By the highest weight vector property we have
\begin{equation} \label{zero}
W(\Lambda_i)=U(\bar{\goth{n}}_{-}) \cdot v_{\Lambda_i}.
\end{equation}
Now we set
\begin{equation}
W(\Lambda_1)' = U(\bar{\goth{n}}_{\leq -2}) \cdot v_{\Lambda_1}.
\end{equation}
Then
\begin{equation} \label{prime}
W(\Lambda_1)' = W(\Lambda_1)
\end{equation}
because $(x_{\alpha} \otimes t^{-1}) \cdot v_{\Lambda_1}=0$ (recall
(\ref{brackets})).

For $i=0,1$ consider the surjective maps
\begin{eqnarray} \label{surj1}
F_{\Lambda_i}: U(\widehat{\goth{g}}) & \longrightarrow & L(\Lambda_i) \\
a &\mapsto& a \cdot v_{\Lambda_i} .\nonumber
\end{eqnarray}
Let us restrict $F_{\Lambda_i}$ to $U(\bar{\goth{n}}_{-})$ and
$F_{\Lambda_1}$ to $U(\bar{\goth{n}}_{\leq -2})$ and denote these
restrictions by $f_{\Lambda_i}$ and $f_{\Lambda_1}'$:
\begin{eqnarray} \label{surj2}
f_{\Lambda_i}: U(\bar{\goth{n}}_{-}) & \longrightarrow & W(\Lambda_i)\\
a & \mapsto & a \cdot v_{\Lambda_i}, \nonumber
\end{eqnarray}
\begin{eqnarray}
f_{\Lambda_1}': U(\bar{\goth{n}}_{\leq -2}) & \longrightarrow &
W(\Lambda_1)' \\ 
a & \mapsto & a \cdot v_{\Lambda_1}. \nonumber
\end{eqnarray}
Our main goal is to describe the kernels $\mbox{Ker} \;
f_{\Lambda_0}$, $\mbox{Ker}\; f_{\Lambda_1}$ and $\mbox{Ker} \;
f_{\Lambda_1}'$. We proceed to do this.

Throughout the rest of the paper we will write $a(m)$ for the action
of $ a \otimes t^m$ on any $\widehat{\goth{g}}$-module, where $a \in
\goth{g}$ and $ m \in \mathbb{Z}$ (recall \cite{CLM1}). In particular, we
have the operator $x_{\alpha}(m)$, the image of $x_{\alpha} \otimes
t^m$. We consider the following formal infinite sums indexed by $t \in
\mathbb{Z}$:
\begin{equation} \label{R_t}
R_t=\sum_{m_1+m_2=-t} x_{\alpha}(m_1)x_{\alpha}(m_2).
\end{equation}
For each $t$, $R_t$ acts naturally on any highest weight
$\widehat{\goth{g}}$-module and, in particular, on $L(\Lambda_i)$,
$i=1,2$.

In order to describe $\mbox{Ker} \; f_{\Lambda_i}$ and $\mbox{Ker} \;
f_{\Lambda_1}'$ we shall truncate each $R_t$ as follows:
\begin{equation} \label{R^0_t}
R_t^0 =\sum_{m_1, m_2 \leq -1, \;
m_1+m_2=-t}x_{\alpha}(m_1)x_{\alpha}(m_2), \; \; t \geq 2.
\end{equation}
We shall often be viewing $R_t^0$ as an element of
$U(\bar{\goth{n}})$, and in fact of $U(\bar{\goth{n}}_{-})$, rather
than as an endomorphism of a module such as $L(\Lambda_i)$. In general
it will be clear from the context when expressions such as
(\ref{R^0_t}) are understood as elements of a universal enveloping
algebra or as operators.  It will also be convenient to take $m_1, m_2
\leq -2$ in (\ref{R_t}), to obtain other elements of
$U(\bar{\goth{n}})$, which we denote by $R_t^1$:
\begin{equation}
R_t^1=\sum_{m_1, m_2 \leq -2, \; m_1+m_2=-t}
x_{\alpha}(m_1)x_{\alpha}(m_2), \; \; t \geq 4.
\end{equation}  

One can view $U(\bar{\goth{n}}_{-})$ and $U(\bar{\goth{n}}_{\leq -2})$
as the polynomial algebras
\begin{equation}
U(\bar{\goth{n}}_{-})= \mathbb{C}[x_{\alpha}(-1), x_{\alpha}(-2), \dots ]
\end{equation}
and
\begin{equation}
U(\bar{\goth{n}}_{\leq-2})= \mathbb{C}[x_{\alpha}(-2), x_{\alpha}(-3), \dots ].
\end{equation}
Then
\begin{equation} \label{deco}
U(\bar{\goth{n}}_{-})= U(\bar{\goth{n}}_{\leq -2}) \oplus
U(\bar{\goth{n}}_{-}) x_{\alpha}(-1)
\end{equation} 
and we have the corresponding projection 
\begin{equation} \label{rho}
\rho: U(\bar{\goth{n}}_{-}) \longrightarrow U(\bar{\goth{n}}_{\leq -2}).
\end{equation}
Notice that
\begin{equation} 
R_t^1= \rho (R_t^0).
\end{equation}

Set
\begin{equation} \label{ideal1}
I_{\Lambda_0}= \sum_{t \geq 2} U(\bar{\goth{n}}_{-})R_t^0 \subset
U(\bar{\goth{n}}_{-})
\end{equation}
and
\begin{equation} \label{ideal2}
I_{\Lambda_1}=\sum_{t \geq 2} U(\bar{\goth{n}}_{-})R_t^0 +
U(\bar{\goth{n}}_{-})x_{\alpha}(-1) \subset U(\bar{\goth{n}}_{-}).
\end{equation}
Observe that (\ref{ideal2}) can be written as 
\begin{equation} \label{ideal22}
I_{\Lambda_1}=\sum_{t \geq 4} U(\bar{\goth{n}}_{-})R_t^1 +
U(\bar{\goth{n}}_{-}) x_{\alpha}(-1).
\end{equation}

\begin{remark}
\em 
Note that
\begin{equation} \label{eq0}
I_{\Lambda_0} \subset I_{\Lambda_1},
\end{equation}
and in fact
\begin{equation}
I_{\Lambda_1}=I_{\Lambda_0} + U(\bar{\goth{n}}_{-}) x_{\alpha}(-1).
\end{equation}

\end{remark}

Also set
\begin{equation} \label{I'}
I_{\Lambda_1}'= \sum_{t \geq 4} U(\bar{\goth{n}}_{\leq -2}) R_t^1
\subset U(\bar{\goth{n}}_{\leq -2}).
\end{equation}
Observe that
\begin{equation}
\rho  (I_{\Lambda_1})= I_{\Lambda_1}',
\end{equation}
and in fact,
\begin{equation} \label{1-1'}
I_{\Lambda_1}= I_{\Lambda_1}' \oplus U(\bar{\goth{n}}_{-})x_{\alpha}(-1).
\end{equation}

It is well known (\cite{B}, \cite{FLM}) that there is a natural vertex
operator algebra structure on $L(\Lambda_0)$ with a vertex operator
map
\begin{eqnarray}
Y(\cdot , x) : L(\Lambda_0) & \longrightarrow & \mbox{End} \; 
L(\Lambda_0) \; [[x, x^{-1}]] \\
v  & \mapsto &Y(v, x)= \sum_{m \in \mathbb{Z}} v_m x^{-m-1}, \nonumber 
\end{eqnarray} 
which satisfies certain conditions, and with $v_{\Lambda_0}$ as vacuum vector. 
In particular, we have
\begin{equation} \label{vertex}
Y(e^{\alpha}, x)= \sum_{m \in \mathbb{Z}} x_{\alpha}(m) x^{-m-1},
\end{equation}
where $e^{\alpha}$ is viewed as an element of $L(\Lambda_0)$.  It is
also well known that $L(\Lambda_1)$ has a natural
$L(\Lambda_0)$-module structure.

The vector spaces $L(\Lambda_i)$ are graded with respect to a standard
action of the Virasoro algebra operator $L(0)$, usually referred to as
grading by {\it conformal weight}. For $m$ an integer,
\begin{equation} \label{wt}
\mbox{wt} \; x_{\alpha}(m)= -m,
\end{equation}
where $x_{\alpha}(m)$ is viewed as either an operator or as an
element of $U(\bar{\goth{n}})$.  For any element $\lambda$ of the
weight lattice of $\goth{g}$ we have
\begin{equation}
\mbox{wt} \; e^{\lambda}= \frac{1}{2} \langle \lambda, \lambda \rangle
\end{equation}
(cf. \cite{CLM1}).
In particular, 
\begin{equation}
\mbox{wt} \; v_{\Lambda_0}=0 \; \; \; \mbox{and} \; \; \; \mbox{wt}
\; v_{\Lambda_1}= \frac{1}{4}
\end{equation}
(recall (\ref{v})).  

The vector spaces $L(\Lambda_i)$ have also a grading given by the
eigenvalues of the operator $\frac{1}{2} \alpha(0)=\frac{1}{2} h(0)$,
called the grading by {\it charge}. This is compatible with the
grading by conformal weight.  For any $m \in \mathbb{Z}$,
$x_{\alpha}(m)$, viewed as either an operator or as an element of
$U(\bar{\goth{n}})$, has charge $1$. Also, $e^{\alpha/2}$, viewed as
either an operator or as an element of $L(\Lambda_1)$, has charge
$\frac{1}{2}$. We shall consider these gradings restricted to the
principal subspaces $W(\Lambda_i)$. For any $m_1, \dots , m_k \in
\mathbb{Z}$,
\begin{equation} \label{v0}
x_{\alpha}(m_1) \cdots x_{\alpha}(m_k) \cdot v_{\Lambda_0} \in W(\Lambda_0)
\end{equation}
and 
\begin{equation} \label{v1}
x_{\alpha}(m_1) \cdots x_{\alpha}(m_k) \cdot v_{\Lambda_1} \in W(\Lambda_1)
\end{equation}
have weights $-m_1-\cdots -m_k$ and $ -m_1-\cdots -m_k+ \frac{1}{4}$,
respectively. Their charges are $k$ and $k+\frac{1}{2}$,
respectively.  See Section 2 of \cite{CLM1} for further details,
background and notation.

\begin{remark} \label{L(0)}
\em
It is clear that
$$
L(0) \; \mbox{\rm Ker} \; f_{\Lambda_i} \subset \mbox{\rm Ker} \;
f_{\Lambda_i}
$$
for $i=0,1$ and
$$
L(0) \; \mbox{Ker} \; f_{\Lambda_1}' \subset \mbox{Ker} \; f_{\Lambda_1}'.
$$
Also, $R_t^0$ and $R_t^1$ have conformal weight $t$:
$$
L(0) R_t^0 = t R_t^0 \; \; \; \mbox{for all} \; \; \; t \geq 2
$$
and
$$
L(0) R_t^1= t R_t^1 \; \; \; \mbox{for all} \; \; \; t \geq 4.
$$
In particular, the subspaces $I_{\Lambda_0}$, $I_{\Lambda_1}$ and
$I_{\Lambda_1}'$ are $L(0)$-stable.  We also have that $R_t^0$ and
$R_t^1$ have charge 2, and the spaces $\mbox{Ker} \; f_{\Lambda_i}$,
$\mbox{Ker} \; f_{\Lambda_1}'$, $I_{\Lambda_i}$ and $I_{\Lambda_1}'$
for $i=0,1$ are graded by charge. Hence these spaces are graded by
both weight and charge, and the two gradings are compatible.
\end{remark}

We will prove the following result describing the kernels of
$f_{\Lambda_i}$ and $f_{\Lambda_1}'$ (recall (\ref{deco}) and
(\ref{1-1'})):
\begin{theorem} \label{th1}
We have
\begin{equation}
\mbox{ \rm Ker} \; f_{\Lambda_0}= I_{\Lambda_0}, \; \; \mbox{ \rm and also,}
\; \; \; \mbox{\rm Ker} \; f_{\Lambda_1}'= I_{\Lambda_1}',
\end{equation}
or equivalently,
\begin{equation}
{\rm Ker} \; f_{\Lambda_1}= I_{\Lambda_1}.
\end{equation}
\end{theorem}

We will actually prove a restatement of this assertion (see Theorem
\ref{th2} below). For this reason we need to introduce the generalized
Verma modules, in the sense of \cite{L1}, \cite{GL}, \cite{L2}, for
$\widehat{\goth{g}}$, as well as what we shall call the principal
subspaces of these generalized Verma modules.

The generalized Verma module $N(\Lambda_0)$ is defined as the induced
$\widehat{\goth{g}}$-module
\begin{equation}
N(\Lambda_0) = U(\widehat{\goth{g}}) \otimes
_{U(\widehat{\goth{g}}_{\geq 0})} \mathbb{C}v_{\Lambda_0}^N,
\end{equation}
where $\goth{g} \otimes \mathbb{C}[t]$ acts trivially and ${\bf k}$
acts as the scalar $1$ on $\mathbb{C}v_{\Lambda_0}^N$;
$v_{\Lambda_0}^N$ is a highest weight vector.  From the
Poincar\'e-Birkhoff-Witt theorem we have
\begin{equation} \label{pbw}
N(\Lambda_0) =U(\widehat{\goth{g}}_{<0}) \otimes_{\mathbb{C}}
U(\widehat{\goth{g}}_{\geq 0}) \otimes _{U(\widehat{\goth{g}}_{\geq
0})} \mathbb{C}v_{\Lambda_0}^N= U(\widehat{\goth{g}}_{<0}) \otimes
_{\mathbb{C}} \mathbb{C}v_{\Lambda_0}^N =U(\widehat{\goth{g}}_{<0}),
\end{equation}
with the natural identifications.  
Similarly define the generalized Verma module
$$
N(\Lambda_1)=U(\widehat{\goth{g}}) \otimes_{U(\widehat{\goth{g}}_{\geq
0})} U,
$$
where $U$ is a two-dimensional irreducible $\goth{g}$-module with a
highest-weight vector $v_{\Lambda_1}^N$, and where $\goth{g} \otimes t
\mathbb{C}[t]$ acts trivially and ${\bf k}$ by $1$.  By the
Poincar\'e-Birkhoff-Witt theorem we have the identification
$$
N(\Lambda_1)= U(\widehat{\goth{g}}_{<0}) \otimes_{\mathbb{C}} U.
$$

For $i=0,1$ we have the natural surjective $\widehat{\goth{g}}$-module maps
\begin{eqnarray} \label{starr}
F_{\Lambda_i}^N: U(\widehat{\goth{g}}) & \longrightarrow & N(\Lambda_i)  \\
a & \mapsto & a \cdot v_{\Lambda_i}^N \nonumber
\end{eqnarray}
(cf. (\ref{surj1})).

\begin{remark} \label{star}
\rm The restriction of (\ref{starr}) to $U(\widehat{\goth{g}}_{<0})$
is a $U(\widehat{\goth{g}}_{<0})$-module isomorphism for $i=0$ and a
$U(\widehat{\goth{g}}_{<0})$-module injection for $i=1$.
\end{remark}

Set  
\begin{equation}
W^N(\Lambda_i)= U(\bar{\goth{n}}) \cdot v_{\Lambda_i}^N,
\end{equation}
an $\bar{\goth{n}}$-submodule of $N(\Lambda_i)$ (cf. (\ref{ps})).  We
shall call $W^N(\Lambda_i)$ the {\it principal subspace} of the
generalized Verma module $N(\Lambda_i)$.  By the highest weight vector
property we have
\begin{equation}
W^N(\Lambda_i) = U(\bar{\goth{n}}_{-}) \cdot v_{\Lambda_i}^N.
\end{equation}
Also consider the subspace  
\begin{equation}
W^N(\Lambda_1)' = U(\bar{\goth{n}}_{\leq -2}) \cdot v_{\Lambda_1}^N
\end{equation}
of $W^N(\Lambda_1)$.

\begin{remark} \label{star-pr}
\rm 
In view of Remark \ref{star}, the maps
\begin{eqnarray} \label{eq1}
U(\bar{\goth{n}}_{-}) & \longrightarrow & W^N(\Lambda_i) \\
a & \mapsto & a \cdot v_{\Lambda_i}^N \nonumber 
\end{eqnarray}
are $\bar{\goth{n}}_{-}$-module isomorphisms and the map
\begin{eqnarray} \label{eq2}
U(\bar{\goth{n}}_{\leq -2}) & \longrightarrow & W^N(\Lambda_1)' \\
a & \mapsto & a \cdot v_{\Lambda_1}^N \nonumber
\end{eqnarray}
is an $\bar{\goth{n}}_{\leq -2}$-module isomorphism.
\end{remark}

In particular, we have the natural identifications
\begin{equation}
W^N(\Lambda_1)' \simeq W^N(\Lambda_1)/ U(\bar{\goth{n}}_{-})
x_{\alpha}(-1) \cdot v_{\Lambda_1}^N \simeq U(\bar{\goth{n}}_{\leq
-2}),
\end{equation}
in view of (\ref{deco}).  

We have the natural surjective $\widehat{\goth{g}}$-module maps
\begin{eqnarray} \label{Pi}
\Pi_{\Lambda_i} : N(\Lambda_i) & \longrightarrow & L(\Lambda_i) \\
a \cdot v_{\Lambda_i}^N & \mapsto & a \cdot v_{\Lambda_i} \nonumber
\end{eqnarray}
for $a \in \widehat{\goth{g}}$, $i=0,1$. Set
\begin{equation} \label{ker}
N^1(\Lambda_i)= \mbox{Ker} \; \Pi_{\Lambda_i}.
\end{equation}
The restrictions of
$\Pi_{\Lambda_i}$ to $W^N(\Lambda_i)$ (respectively,
$W^N(\Lambda_1)'$) are $\bar{\goth{n}}$-module (respectively,
$\bar{\goth{n}}_{\leq -2}$-module) surjections:
\begin{equation} \label{pi}
\pi_{\Lambda_i}: W^N(\Lambda_i) \longrightarrow W(\Lambda_i)
\end{equation}
for $i=0,1$ and 
\begin{equation}
\pi_{\Lambda_1}': W^N(\Lambda_1)' \longrightarrow W(\Lambda_1)
\end{equation}
(recall (\ref{prime})).

As in the case of $L(\Lambda_i)$, there are natural actions of the
Virasoro algebra operator $L(0)$ on $N(\Lambda_i)$ for $i=0,1$, giving
gradings by conformal weight, and these spaces are also compatibly
graded by charge, by means of the operator $\frac{1}{2}
\alpha(0)=\frac{1}{2}h(0)$.  We shall restrict these gradings to the
principal subspaces $W^N(\Lambda_i)$. The elements of $W^N(\Lambda_i)$
given by (\ref{v0}) and (\ref{v1}) with $v_{\Lambda_i}$ replaced by
$v_{\Lambda_i}^N$ have the same weights and charges as in those cases.

\begin{remark} \label{L(0)k} 
\rm Since the maps $\pi_{\Lambda_i}$, $i=0,1$, and $\pi_{\Lambda_1}'$
commute with the actions of $L(0)$, the kernels $\mbox{Ker} \;
\pi_{\Lambda_i}$ and $\mbox{Ker} \; \pi_{\Lambda_1}'$ are
$L(0)$-stable.  These maps preserve charge as well, so that
$\mbox{Ker} \; \pi_{\Lambda_i}$ and $\mbox{Ker} \; \pi_{\Lambda_1}'$
are also graded by charge.
\end{remark}

In view of Remark \ref{star}, the following assertion is equivalent to that of
Theorem \ref{th1}:
\begin{theorem} \label{th2}
We have
\begin{equation} \label{eqone} 
{\rm Ker} \; \pi_{\Lambda_0} = I_{\Lambda_0} \cdot v_{\Lambda_0}^N \;
\; (\subset N^1(\Lambda_0)).
\end{equation}
Moreover, 
\begin{equation} \label{eqtwo}
{\rm Ker} \; \pi_{\Lambda_1}' = I_{\Lambda_1}' \cdot v_{\Lambda_1}^N
\; \; (\subset N^1(\Lambda_1)),
\end{equation}
or equivalently,
\begin{equation} \label{eqthree}
{\rm Ker} \; \pi_{\Lambda_1}= I_{\Lambda_1} \cdot v_{\Lambda_1}^N \;
\; (\subset N^1(\Lambda_1)).
\end{equation}
\end{theorem}

\section{ Proof of the main result}
\setcounter{equation}{0}
Continuing to use the setting of \cite{CLM1}, we have
\[
V_P=L(\Lambda_0) \oplus L(\Lambda_1).
\]
Consider the linear isomorphism 
\begin{equation} \label{map}
e^{\alpha/2} : V_P \longrightarrow V_P.
\end{equation}
Its restriction to the principal subspace $W(\Lambda_0)$ of $L(\Lambda_0)$ is
\begin{equation} \label{shift}
e^{\alpha/2} : W(\Lambda_0)  \longrightarrow W(\Lambda_1).
\end{equation}
We have
\begin{equation} \label{op}
e^{\alpha/2} \; x_{\alpha}(m)=x_{\alpha}(m-1) \; e^{\alpha/2}
\end{equation}
on $V_P$ for $m \in \mathbb{Z}$, and
\begin{equation} \label{vec}
e^{\alpha/2} \cdot v_{\Lambda_0}= v_{\Lambda_1},
\end{equation}
so that
\begin{equation} \label{def-shift}
e^{\alpha/2} \; (x_{\alpha}(m_1)\cdots x_{\alpha}(m_k) \cdot
v_{\Lambda_0})= x_{\alpha}(m_1-1) \cdots x_{\alpha}(m_k-1) \cdot
v_{\Lambda_1}
\end{equation}
for any $m_1, \dots, m_k \in \mathbb{Z}$. Then in particular,
(\ref{shift}) is a linear isomorphism.

We shall follow \cite{CLM1} to construct a lifting 
\begin{equation} \label{hatshift}
\widehat{e^{\alpha/2}}: W^N(\Lambda_0) \longrightarrow W^N(\Lambda_1)'
\end{equation}
of 
$$
e^{\alpha/2}: W(\Lambda_0) \longrightarrow W(\Lambda_1)
$$
in the sense that the following diagram will commute:
$$
\CD
W^N(\Lambda_0) @> \widehat{e^{\alpha/2}}> \sim > W^N(\Lambda_1)'  \\
@V\pi_{\Lambda_0}VV @V\pi_{\Lambda_1}'VV  \\
W(\Lambda_0)\ @> {e^{\alpha/2}}>\sim >
W(\Lambda_1).
\endCD
$$ 
However, here we shall use our current notation, which involves
generalized Verma modules and principal subspaces of these modules.

For any integers $m_1, \dots, m_k <0$ we set
\begin{equation} \label{def-hatshift}
\widehat{e^{\alpha/2}} \; ( x_{\alpha}(m_1) \cdots x_{\alpha}(m_k)
\cdot v_{\Lambda_0}^N)= x_{\alpha}(m_1-1) \cdots x_{\alpha}(m_k-1)
\cdot v_{\Lambda_1}^N,
\end{equation}
which is well defined, since $U(\bar{\goth{n}}_{-})$, viewed as the
polynomial algebra
$$\mathbb{C}[x_{\alpha}(-1), x_{\alpha}(-2), \dots],$$
maps isomorphically onto $W^N(\Lambda_0)$ under the map
(\ref{eq1}). Since
$$
W^N(\Lambda_1)'=U(\bar{\goth{n}}_{\leq -2}) \cdot v_{\Lambda_1}^N=
\widehat{e^{\alpha/2}} \; (U(\bar{\goth{n}}_{-}) \cdot
v_{\Lambda_0}^N),
$$
the linear map $\widehat{e^{\alpha/2}}$ is surjective. By Remark
\ref{star-pr} it is also injective. Denote by
\begin{equation} \label{shift1}
(\widehat{e^{\alpha/2}})^{-1}=\widehat{e^{-\alpha/2}}: W^N(\Lambda_1)'
\longrightarrow W^N(\Lambda_0)
\end{equation}
its inverse. Then $\widehat{e^{\alpha/2}}$ is indeed a lifting of
$e^{\alpha/2}$ as desired, and $\widehat{e^{-\alpha/2}}$ is
correspondingly a lifting of the inverse
\begin{equation} \label{shift11}
e^{-\alpha/2}: W(\Lambda_1) \longrightarrow W(\Lambda_0).
\end{equation}

\begin{lemma} \label{lemma1} We have
\begin{equation} \label{0=1}
\widehat{e^{\alpha/2}} \; (I_{\Lambda_0} \cdot v_{\Lambda_0}^N)=
I_{\Lambda_1}' \cdot v_{\Lambda_1}^N.
\end{equation}
\end{lemma} 
{\bf Proof:} 
For any $t \geq 2$,
\begin{eqnarray} 
\widehat{e^{\alpha/2}} \; (R_t^0 \cdot v_{\Lambda_0}^N) & =&
\widehat{e^{\alpha/2}} \left ( \sum_{m_1, m_2 \leq -1, \; m_1+m_2=-t}
x_{\alpha}(m_1) x_{\alpha}(m_2) \cdot v_{\Lambda_0}^N \right )
\nonumber \\
& =& \sum_{m_1, m_2 \leq -1, \; m_1+m_2=-t} x_{\alpha}(m_1-1)
x_{\alpha}(m_2-1) \cdot v_{\Lambda_1}^N= R_{t+2}^1 \cdot
v_{\Lambda_1}^N. \nonumber
\end{eqnarray}
In view of the definition of the map $\widehat{e^{\alpha/2}}$ and of
the descriptions (\ref{ideal1}) and (\ref{I'}) of the ideals
$I_{\Lambda_0}$ and $I_{\Lambda_1}'$ we have (\ref{0=1}). $\;\;\; \Box$
\vspace{1em}

Now we restrict (\ref{map}) to the principal subspace $W(\Lambda_1)$
and obtain
\begin{equation} \label{shift2}
e^{\alpha/2} : W(\Lambda_1) \longrightarrow W(\Lambda_0),
\end{equation}
which is only an injection. Since
\begin{equation} \label{e^}
e^{\alpha/2} \cdot v_{\Lambda_1}= e^{\alpha} \cdot v_{\Lambda_0}=
x_{\alpha}(-1) \cdot v_{\Lambda_0}
\end{equation}
(recall (\ref{vec}) and \cite{CLM1}), we have, using (\ref{op}),
\begin{equation} \label{def-shift2}
e^{\alpha/2} \; (x_{\alpha}(m_1) \cdots x_{\alpha}(m_k) \cdot
v_{\Lambda_1}) = x_{\alpha}(m_1-1) \cdots x_{\alpha}(m_k-1)
x_{\alpha}(-1) \cdot v_{\Lambda_0}
\end{equation}
for any $m_1, \dots , m_k \in \mathbb{Z}$. As above, we construct a
natural lifting
\begin{equation} \label{shift22}
\widehat{e^{\alpha/2}}: W^N(\Lambda_1)'  \longrightarrow W^N(\Lambda_0)
\end{equation}
of the map (\ref{shift2}), making the diagram
$$
\CD
W^N(\Lambda_1)' @> \widehat{e^{\alpha/2}} >> W^N(\Lambda_0) \\
@V\pi_{\Lambda_1}'VV @V\pi_{\Lambda_0}VV  \\
W(\Lambda_1)\ @> {e^{\alpha/2}}> >
W(\Lambda_0)
\endCD
$$ 
commute, 
by taking 
\begin{equation} \label{def-shift22}
\widehat{e^{\alpha/2}} \; (x_{\alpha}(m_1) \cdots x_{\alpha}(m_k)
\cdot v_{\Lambda_1}^N)= x_{\alpha}(m_1-1) \cdots
x_{\alpha}(m_k-1)x_{\alpha}(-1) \cdot v_{\Lambda_0}^N
\end{equation}
for $m_1, \dots, m_k \leq -2$.  Then, like (\ref{shift2}), the map
(\ref{shift22}) is an injection and not a surjection.

\begin{lemma}\label{lemma2} We have
\begin{equation}
\widehat{e^{\alpha/2}} \; (I_{\Lambda_1}' \cdot v_{\Lambda_1}^N)
\subset I_{\Lambda_0} \cdot v_{\Lambda_0}^N.
\end{equation} 
\end{lemma}
{\bf Proof:} 
As in the proof of the previous lemma, we calculate 
that for any integer $t \geq 4$,
$$
\widehat{e^{\alpha/2}} \; (R_t^1 \cdot v_{\Lambda_1}^N)= R_{t+2}^0
x_{\alpha}(-1)\cdot v_{\Lambda_0}^N - x_{\alpha}(-t) R_3^0 \cdot
v_{\Lambda_0}^N- 2x_{\alpha}(-t-1)R_2^0 \cdot v_{\Lambda_0}^N,
$$
which proves our lemma.
$\;\;\; \Box$
\vspace{1em}

Our next goal is to prove the main result, Theorem \ref{th1}, or
equivalently, Theorem \ref{th2}. It is sufficient to prove
(\ref{eqone}) and (\ref{eqtwo}).

We notice first that
\begin{equation} \label{first inclusions}
I_{\Lambda_0} \cdot v_{\Lambda_0}^N \subset \mbox{Ker} \;
\pi_{\Lambda_0} \; \; \mbox{and} \; \; I_{\Lambda_1}' \cdot
v_{\Lambda_1}^N \subset \mbox{Ker} \; \pi_{\Lambda_1}'.
\end{equation}
Indeed, as is well known, the square of the vertex operator
$Y(e^{\alpha}, x)$ is well defined (the components $x_{\alpha}(m)$, $m
\in \mathbb{Z}$, of this vertex operator commute) and equals zero on
$L(\Lambda_i)$, and in particular on $W(\Lambda_i)$, for $i=1,2$. The
expansion coefficients of $Y(e^{\alpha}, x)^2$ are the operators
$R_{-t}$, $t \in \mathbb{Z}$:
\begin{equation}
Y(e^{\alpha}, x)^2= \sum_{t \in \mathbb{Z}} \left ( \sum_{m_1+m_2=t}
x_{\alpha}(m_1)x_{\alpha}(m_2)\right )x^{-t-2}
\end{equation}
(recall (\ref{R_t}) and (\ref{vertex})), and this proves (\ref{first
inclusions}).

Now we define a shift, or translation, automorphism
\begin{equation} \label{translation}
\tau: U(\bar{\goth{n}}) \longrightarrow U(\bar{\goth{n}})
\end{equation}
by
$$
\tau (x_{\alpha}(m_1) \cdots x_{\alpha}(m_k))= x_{\alpha}(m_1-1)
\cdots x_{\alpha}(m_k-1)
$$
for any integers $m_1, \dots, m_k$; since $U(\bar{\goth{n}}) \simeq
\mathbb{C}[x_{\alpha}(m), \; m \in \mathbb{Z}]$, the map
(\ref{translation}) is well defined. For any integer $s$, the
$s^{\mbox{th}}$ power
\begin{equation}
\tau ^s: U(\bar{\goth{n}}) \longrightarrow U(\bar{\goth{n}})
\end{equation}
is given by
$$
\tau ^s (x_{\alpha}(m_1) \cdots x_{\alpha}(m_k))= x_{\alpha}(m_1-s)
\cdots x_{\alpha}(m_k-s).
$$

\begin{remark}
\rm Assume that $a \in U(\bar{\goth{n}})$ is a nonzero element
homogeneous with respect to both the weight and charge gradings. We
also assume that $a$ has positive charge, so that $a$ is a linear
combination of monomials of a fixed positive charge (degree). Then
$\tau^s(a)$ has the same properties, and
\begin{equation} \label{greater}
\mbox{wt} \; \tau^s (a) > \mbox{wt} \; a \; \; \; \mbox{for} \; \; \; s > 0
\end{equation}
and
\begin{equation} \label{smaller}
\mbox{wt} \; \tau^s (a) < \mbox{wt} \; a \; \; \; \mbox{for} \; \; \; s <0,
\end{equation}
by (\ref{wt}). If $a$ is a constant, that is, $a$ has charge zero,
then of course
\begin{equation} \label{constant}
\tau^s(a)= a
\end{equation} 
and $\tau^s(a)$ and $a$ have the same weight  and charge.
\end{remark}

We also have:

\begin{remark} \label{rem} 
\rm Using the map $\tau$ we can re-express (\ref{shift}),
(\ref{hatshift}), (\ref{shift2}) and (\ref{shift22}) as follows:
\begin{equation} \label{t1}
e^{\alpha/2} (a \cdot v_{\Lambda_0})= \tau(a) \cdot v_{\Lambda_1}, \;
\; a \in U(\bar{\goth{n}}),
\end{equation}
\begin{equation} \label{t2}
\widehat{e^{\alpha/2}} ( a \cdot v_{\Lambda_0}^N) = \tau (a) \cdot
v_{\Lambda_1}^N, \; \; a \in U(\bar{\goth{n}}_{-}),
\end{equation}
\begin{equation} \label{t3}
e^{\alpha/2} (a \cdot v_{\Lambda_1})= \tau(a) x_{\alpha}(-1) \cdot
v_{\Lambda_0}, \; \; a \in U(\bar{\goth{n}})
\end{equation}
and
\begin{equation} \label{t4}
\widehat{e^{\alpha/2}} ( a \cdot v_{\Lambda_1}^N)= \tau(a)
x_{\alpha}(-1) \cdot v_{\Lambda_0}^N, \; \; a \in
U(\bar{\goth{n}}_{\leq -2})
\end{equation}
(recall (\ref{def-shift}), (\ref{def-hatshift}), (\ref{def-shift2})
and (\ref{def-shift22})). The corresponding inverse maps of course
involve $\tau^{-1}$.
\end{remark}

\begin{lemma} \label{small}
We have
$$
\tau ( I_{\Lambda_0}) \subset I_{\Lambda_0} + U(\bar{\goth{n}}_{-})
x_{\alpha}(-1)=I_{\Lambda_1}.
$$
\end{lemma}
{\bf Proof:} Let $t \geq 2$. Then
$$
\tau (R_t^0)= R_{t+2}^0 - 2x_{\alpha}(-t-1) x_{\alpha}(-1). \; \; \; \; \Box
$$

Recall from \cite{CLM1} the role of intertwining vertex operators among
triples of the $L(\Lambda_0)$-modules $L(\Lambda_0)$ and
$L(\Lambda_1)$, and the role of certain terms and factors of these
intertwining operators. We will use the constant term of the
intertwining operator ${\cal Y}(e^{\alpha/2}, x)$ of type
$$ \left(
\begin{array} {c}
L( \Lambda_1)                    \\
\begin{array}{cc}
L(\Lambda_1)   &  L(\Lambda_0)
\end{array} 
\end{array}
\right), $$
and denote it by ${\cal Y}_c(e^{\alpha/2}, x)$. This is the operator
denoted by $o(e^{\alpha/2})$ in Theorem 4.1 of \cite{CLM1}. Then
\begin{equation} \label{Y_c}
{\cal Y}_c(e^{\alpha/2}, x): W(\Lambda_0) \longrightarrow
W(\Lambda_1),
\end{equation}
it sends $v_{\Lambda_0}$ to $v_{\Lambda_1}$, and it commutes with the
action of $\bar{\goth{n}}$ (see the beginning of the proof of Theorem
4.1 in \cite{CLM1}).

{\bf Proof of Theorem \ref{th2}:} 
In view of (\ref{first inclusions}), it is sufficient to prove 
\begin{equation}
\mbox{Ker} \; \pi_{\Lambda_0} \subset I_{\Lambda_0} \cdot
v_{\Lambda_0}^N \; \; \mbox{and} \; \; \mbox{Ker} \; \pi_{\Lambda_1}'
\subset I_{\Lambda_1}' \cdot v_{\Lambda_1}^N.
\end{equation}

First we show that (\ref{eqtwo}) follows from (\ref{eqone}), whose
truth we now assume.  Let $a \cdot v_{\Lambda_1}^N \in \mbox{Ker} \;
\pi_{\Lambda_1}'$, where $a \in U(\bar{\goth{n}}_{\leq -2})$, so that
$a \cdot v_{\Lambda_1} =0$ in $W(\Lambda_1)$. By (\ref{t1}) we have
$\tau^{-1}(a) \cdot v_{\Lambda_0} =0$ in $W(\Lambda_0)$, and so
$$
\tau^{-1}(a) \cdot v_{\Lambda_0}^N \in \mbox{Ker} \;
\pi_{\Lambda_0}=I_{\Lambda_0} \cdot v_{\Lambda_0}^N
$$
by assumption.
Thus
$$
\widehat{e^{\alpha/2}} \; (\tau^{-1}(a) \cdot v_{\Lambda_0}^N) \in
\widehat{e^{\alpha/2}} \; (I_{\Lambda_0} \cdot v_{\Lambda_0}^N),$$ and
by (\ref{t2}) and Lemma \ref{lemma1} we obtain
$$
a \cdot v_{\Lambda_1}^N \in I_{\Lambda_1}' \cdot v_{\Lambda_1}^N.
$$
Thus we have the inclusion $\mbox{Ker} \; \pi_{\Lambda_1}' \subset
I_{\Lambda_1}' \cdot v_{\Lambda_1}^N$, and so we have (\ref{eqtwo}).

It remains to prove
\begin{equation} \label{inclusion}
\mbox{Ker} \; \pi_{\Lambda_0} \subset I_{\Lambda_0} \cdot v_{\Lambda_0}^N.
\end{equation}
We will prove by contradiction that any element of
$U(\bar{\goth{n}}_{-}) \cdot v_{\Lambda_0}^N$ that lies in
$\mbox{Ker} \; \pi_{\Lambda_0}$ lies in $I_{\Lambda_0} \cdot
v_{\Lambda_0}^N$ as well.  Suppose then that there exists $a \in
U(\bar{\goth{n}}_{-})$ such that
\begin{equation} \label{contra}
a \cdot v_{\Lambda_0}^N \in \mbox{Ker} \; \pi_{\Lambda_0} \; \; \;
\mbox{but} \; \; \; a \cdot v_{\Lambda_0}^N \notin I_{\Lambda_0}
\cdot v_{\Lambda_0}^N.
\end{equation}
By Remarks \ref{L(0)} and \ref{L(0)k} we may and do assume that $a$
is homogeneous with respect to the weight and charge gradings. Note
that $a$ is nonzero and in fact nonconstant, and so it has positive
weight and charge.  We choose such an element $a$ of smallest possible
weight that satisfies (\ref{contra}).

We claim that there is in fact an element lying in
$U(\bar{\goth{n}}_{-}) x_{\alpha}(-1)$ and in addition having all of
the properties of $a$.

By (\ref{deco}) we have a unique decomposition $a= b x_{\alpha}(-1)
+c$ with $b \in U(\bar{\goth{n}}_{-})$ and $c \in
U(\bar{\goth{n}}_{\leq -2})$. The elements $b$ and $c$ are homogeneous
with respect to the weight and charge gradings. In fact,
\begin{equation} \label{weights}
\mbox{wt} \; b = \mbox{wt} \; a -1, \; \; \; \mbox{wt} \; c = \mbox{wt} \; a,
\end{equation}
and similarly, the charge of $b$ is one less than that of $a$ and the
charges of $c$ and $a$ are equal.  We have
\begin{equation} \label{tau inverse a}
\tau^{-1} (a)= \tau^{-1} (b) x_{\alpha}(0)+ \tau^{-1} (c)
\end{equation}
with 
\begin{equation} \label{tau inverse c}
\tau^{-1} (c) \in U(\bar{\goth{n}}_{-}).
\end{equation}

Since $a \in U(\bar{\goth{n}}_{-})$ and 
$$
a \cdot v_{\Lambda_0}=0,
$$
by applying the linear map ${\mathcal Y}_c(e^{\alpha/2}, x)$ and using
its properties (recall (\ref{Y_c})) we have
\begin{equation} \label{v_1}
a \cdot v_{\Lambda_1}=0;
\end{equation}
what we have just observed is that
\begin{equation}
\mbox{Ker} \; f_{\Lambda_0} \subset \mbox{Ker} \; f_{\Lambda_1}
\end{equation}
(recall the notation (\ref{surj2})).
By (\ref{v_1}) and (\ref{t1}) we obtain
$$
\tau^{-1} (a) \cdot v_{\Lambda_0}=0 \; \; \mbox{in} \; \; W(\Lambda_0),
$$
and so
\begin{equation} \label{gen1}
\tau^{-1}(a) \cdot v_{\Lambda_0}^N \in \mbox{Ker} \; \pi_{\Lambda_0};
\end{equation}
this shows that if $u$ lies in $U(\bar{\goth{n}})$ (not necessarily in
$U(\bar{\goth{n}}_{-})$), then
\begin{equation}
u \cdot v_{\Lambda_1}^N \in \mbox{Ker} \; \pi_{\Lambda_1} \Rightarrow
\tau^{-1} (u) \cdot v_{\Lambda_0}^N \in \mbox{Ker} \; \pi_{\Lambda_0}.
\end{equation}
But (\ref{tau inverse a}) and (\ref{gen1}) imply
\begin{equation} \label{gen1 c}
\tau^{-1} (c) \cdot v_{\Lambda_0}^N \in \mbox{Ker} \; \pi_{\Lambda_0}.
\end{equation}

We now show that
\begin{equation} \label{gen c}
\tau^{-1} (c) \cdot v_{\Lambda_0}^N   \in I_{\Lambda_0} \cdot v_{\Lambda_0}^N.
\end{equation}
If (\ref{gen c}) did not hold, then since $\tau^{-1} (c) \in
U(\bar{\goth{n}}_{-})$ is doubly homogeneous and
$$
\mbox{wt} \; \tau^{-1} (c) < \mbox{wt} \; c= \mbox{wt} \; a
$$ 
(recall (\ref{smaller}) and (\ref{weights})), we would be
contradicting our assumption that $a \in U(\bar{\goth{n}}_{-})$ is
doubly homogeneous of smallest possible weight satisfying
(\ref{contra}).  Hence (\ref{gen c}) holds, and from (\ref{tau inverse
c}) and Remark \ref{star-pr} we obtain
$$
\tau^{-1} (c) \in I_{\Lambda_0}.
$$ 
Now Lemma \ref{small} applied to $\tau^{-1}(c)$ implies
\begin{equation}
c \in I_{\Lambda_0} +U(\bar{\goth{n}}_{-}) x_{\alpha}(-1).
\end{equation}

In view of Remark \ref{L(0)} we have $c =d+e$, where $d \in
I_{\Lambda_0}$ and $e \in U(\bar{\goth{n}}_{-}) x_{\alpha}(-1)$ are
homogeneous with respect to both gradings and in fact $\mbox{wt} \; d=
\mbox{wt} \; e =\mbox{wt} \; c= \mbox{wt} \; a$, and similarly for
charge. Now $ a= (bx_{\alpha}(-1)+e)+d$, and since $ d \cdot
v_{\Lambda_0}^N \in I_{\Lambda_0} \cdot v_{\Lambda_0}^N \subset
\mbox{Ker} \; \pi_{\Lambda_0}$ and $a \cdot v_{\Lambda_0}^N$ satisfies
(\ref{contra}), we have
\begin{equation} \label{d}
(bx_{\alpha}(-1)+e) \cdot v_{\Lambda_0}^N \in \mbox{Ker} \;
\pi_{\Lambda_0} \; \; \; \mbox{but} \; \; \; (b x_{\alpha}(-1) +e)
\cdot v_{\Lambda_0}^N \notin I_{\Lambda_0} \cdot v_{\Lambda_0}^N.
\end{equation}
Thus we have proved our claim, by means of the doubly homogeneous
element $bx_{\alpha}(-1)+e \in U(\bar{\goth{n}}_{-})x_{\alpha}(-1)$,
which has the same weight as $a$.

In view of the claim, all we have left to do is to show that there
cannot exist an element $a \in U(\bar{\goth{n}}_{-})x_{\alpha}(-1)$
satisfying (\ref{contra}), homogeneous with respect to both the weight
and charge gradings, and of smallest possible weight among all the
elements of $U(\bar{\goth{n}}_{-})$ satisfying (\ref{contra}).

Take such an element $a$. We may in fact assume that
$$a=bx_{\alpha}(-1) \; \; \mbox{with} \; \; b \in
U(\bar{\goth{n}}_{\leq-3}),
$$
where $\bar{\goth{n}}_{\leq -3}$ is the Lie subalgebra of
$\bar{\goth{n}}_{-}$ defined in the obvious way. Indeed, we have the
decomposition
$$
U(\bar{\goth{n}}_{-}) =U(\bar{\goth{n}}_{\leq -3}) \oplus
(U(\bar{\goth{n}}_{-}) x_{\alpha}(-2)+ U(\bar{\goth{n}}_{-})
x_{\alpha}(-1))
$$
into doubly-graded subspaces, and both $x_{\alpha}(-1)^2 \cdot
v_{\Lambda_0}^N$ and $x_{\alpha}(-2)x_{\alpha}(-1) \cdot
v_{\Lambda_0}^N$ are in $I_{\Lambda_0} \cdot v_{\Lambda_0}^N$.

Note that since $b \in U(\bar{\goth{n}}_{\leq -3})$, 
\begin{equation} \label{tau inverse b}
\tau^{-1} (b) \in U(\bar{\goth{n}}_{\leq -2}) \; \; \; \mbox{and} \;
\; \; \tau^{-2} (b) \in U(\bar{\goth{n}}_{-}).
\end{equation}

We also note that
\begin{equation} \label{weight b}
\mbox{wt} \; b =\mbox{wt} \; a-1
\end{equation}
and that the charge of $b$ is one less than that of $a$.

Since $a \cdot v_{\Lambda_0}^N = b x_{\alpha}(-1) \cdot
v_{\Lambda_0}^N \in \mbox{Ker} \; \pi_{\Lambda_0}$, we have $b
x_{\alpha}(-1) \cdot v_{\Lambda_0} =0$ in $W(\Lambda_0)$.  Then by the
second equality in (\ref{e^}) we have
\[
b e^{\alpha} \cdot v_{\Lambda_0}=0 \; \; \mbox{in} \; \; W(\Lambda_0).
\]
By using (\ref{op}) twice we obtain
$$
\tau^{-2}(b) \cdot v_{\Lambda_0} =0 \; \; \mbox{in} \; \; W(\Lambda_0),
$$
so that
\begin{equation} \label{one}
\tau^{-2}(b) \cdot v_{\Lambda_0}^N \in \mbox{Ker} \; \pi_{\Lambda_0}
\; \; \mbox{in} \; \; W^N(\Lambda_0).
\end{equation}

On the other hand,
\begin{equation} \label{two}
\tau^{-2}(b) \cdot v_{\Lambda_0}^N \notin I_{\Lambda_0} \cdot v_{\Lambda_0}^N.
\end{equation}
Indeed, if $\tau^{-2} ( b) \cdot v_{\Lambda_0}^N \in I_{\Lambda_0}
\cdot v_{\Lambda_0}^N$, then
$$\widehat{e^{\alpha/2}}\; ( \widehat{e^{\alpha/2}}\; ( \tau^{-2}(b)
\cdot v_{\Lambda_0}^N)) \in \widehat{e^{\alpha/2}} \;
(\widehat{e^{\alpha/2}} \; (I_{\Lambda_0} \cdot v_{\Lambda_0}^N))
\subset I_{\Lambda_0} \cdot v_{\Lambda_0}^N$$ by Lemmas \ref{lemma1}
and \ref{lemma2}. However, by (\ref {t2}), (\ref{t4}) and (\ref{tau
inverse b}), we have
$$
\widehat{e^{\alpha/2}} \; (\widehat{e^{\alpha/2}} \; (\tau^{-2}(b)
\cdot v_{\Lambda_0}^N))= \widehat{e^{\alpha/2}} \; (\tau^{-1}(b) \cdot
v_{\Lambda_1}^N)= b x_{\alpha}(-1) \cdot v_{\Lambda_0}^N=a \cdot
v_{\Lambda_0}^N.
$$
Thus we obtain $a \cdot v_{\Lambda_0}^N \in I_{\Lambda_0} \cdot
v_{\Lambda_0}^N$, contradicting (\ref{contra}). This proves
(\ref{two}).

But
\begin{equation} \label{three}
\mbox{wt} \; \tau^{-2}(b) \leq \mbox{wt} \; b = \mbox{wt} \; a-1<
\mbox{wt} \; a
\end{equation}
(recall (\ref{smaller}) in case $b$ has positive charge and
(\ref{constant}) when $b$ is a constant, and also (\ref{weight b})).
Thus we have constructed a doubly homogeneous element, namely,
$\tau^{-2} (b) \in U(\bar{\goth{n}}_{-})$, of weight less than that of
$a$ and satisfying (\ref{contra}). This contradiction proves our
theorem and hence Theorem \ref{th1} as well. $\;\;\; \Box$

\begin{remark}
\rm As an immediate consequence of Theorem \ref{th1}, we observe that
any nonzero doubly homogeneous element $a \in
U(\bar{\goth{n}}_{-})$ such that $a \in \mbox{Ker} \; f_{\Lambda_0}=
I_{\Lambda_0} $ has charge at least $2$; that is, there are no nonzero
linear elements $x_{\alpha}(m)$ in $\mbox{Ker} \; f_{\Lambda_0}$, $m
\leq -1$.  Similarly, there are no nonzero linear elements
$x_{\alpha}(m)$ in $\mbox{Ker} \; f_{\Lambda_1}'$ for $m \leq -2$.  We
observe similarly that any homogeneous element of charge $2$ that lies
in $\mbox{Ker} \; f_{\Lambda_0}$ (respectively, $\mbox{Ker} \;
f_{\Lambda_1}'$) is a multiple of $R_t^0$ for some $t \geq 2$
(respectively, $R_t^1$ for some $t \geq 4$).
\end{remark}

\section{ A further reformulation}
\setcounter{equation}{0} 

In this section we shall present a reformulation of formula
(\ref{eqone}) of Theorem \ref{th2} in terms of ideals of vertex
algebras. We shall refer to \cite{LL} for definitions and results
regarding ideals of vertex (operator) algebras and regarding vertex
operator algebra and module structure on generalized Verma modules for
$\widehat{\goth{g}}$.

It is known that $N(\Lambda_0)$ carries a natural structure of vertex
operator algebra, with a vertex operator map
\begin{eqnarray} \nonumber
Y(\cdot, x): N(\Lambda_0) & \longrightarrow & \mbox{End} \;
N(\Lambda_0) [[x, x^{-1}]] \nonumber \\
v & \mapsto & Y(v, x)= \sum_{m \in \mathbb{Z}} v_{m} x^{-m-1} \nonumber
\end{eqnarray}
satisfying certain conditions, with $v_{\Lambda_0}^N$ as vacuum
vector, and with a suitable conformal vector (see Theorem 6.2.18 in
\cite{LL}). The conformal vector gives rise to the Virasoro algebra
operators $L(m)$, $m \in \mathbb{Z}$, in the standard way. Recall
$L(0)$ from Section 2. Here we use the Virasoro algebra operator
$L(-1)$, which acts as a linear operator on $N(\Lambda_0)$ such that
\begin{equation} \label{L(-1)}
L(-1) v_{\Lambda_0}^N= 0 
\end{equation}
and 
\begin{equation} \label{bracket}
[L(-1), Y(v, x)]= \frac{d}{dx} Y(v, x)
\end{equation}
(cf. \cite{LL}). Formula (\ref{bracket}) gives
\begin{equation} \label{b}
[L(-1), v_m]=-mv_{m-1} \; \; \; \mbox{for} \; \; \; v \in
N(\Lambda_0), \; m \in \mathbb{Z}.
\end{equation}
The vector space $N(\Lambda_1)$ has a natural module structure for the
vertex operator algebra $N(\Lambda_0)$, as described in Theorem 6.2.21
of \cite{LL}.  Formulas (\ref{bracket}) and (\ref{b}) hold on
$N(\Lambda_1)$ as well.  Standard arguments show that $W^N(\Lambda_0)$
is a vertex subalgebra of $N(\Lambda_0)$ (cf. the last part of the
proof of Proposition \ref{4} below); $W^N(\Lambda_0)$ does not contain
the conformal vector of $N(\Lambda_0)$. Moreover, $W^N(\Lambda_1)$ is
a $W^N(\Lambda_0)$-submodule of $N(\Lambda_1)$. Also, $L(0)$ preserves
$W^N(\Lambda_i)$ for $i=0,1$, and $L(-1)$ preserves $W^N(\Lambda_0)$.

Recall from (\ref{Pi}) and (\ref{ker}) the natural surjective
$\widehat{\goth{g}}$-module maps $\Pi_{\Lambda_i}$, $i=0,1$, and their
kernels $N^1(\Lambda_i)$.  Then $N^1(\Lambda_i)$ is the unique maximal
proper ($L(0)$-graded) $\widehat{\goth{g}}$-submodule of
$N(\Lambda_i)$ and
\begin{equation} \label{lie}
N^1(\Lambda_0)= U(\widehat{\goth{g}}) x_{\alpha}(-1)^2\cdot
v_{\Lambda_0}^N= U(\mathbb{C}x_{-\alpha} \oplus \goth{g} \otimes
t^{-1} \mathbb{C}[t^{-1}]) x_{\alpha}(-1)^2 \cdot v_{\Lambda_0}^N,
\end{equation}
\begin{equation}
N^1(\Lambda_1)= U(\widehat{\goth{g}}) x_{\alpha}(-1) \cdot
v_{\Lambda_1}^N= U(\mathbb{C}x_{-\alpha} \oplus \goth{g} \otimes
t^{-1} \mathbb{C}[t^{-1}])x_{\alpha}(-1) \cdot v_{\Lambda_1}^N
\end{equation}
(cf. \cite{K}, \cite{LL}).

For the reader's convenience, we recall the definition of the notion
of ideal of a vertex algebra (Definition 3.9.7 in \cite{LL}):

\begin{definition} \label{id}
\rm
An ideal of the vertex algebra $V$ is a subspace $I$ such that for all
$v \in V$ and $w \in I$,
\begin{equation} \label{co1}
Y(v,x)w \in I((x))
\end{equation}
and
\begin{equation} \label{co2}
Y(w, x)v \in I((x)),
\end{equation}
that is, $v_n w \in I$ and $w_n v \in I$ for all $v \in V$, $w \in I$
and $n \in \mathbb{Z}$.
\end{definition}

\begin{remark} \label{skew}
\rm
In view of the skew-symmetry property
$$
Y(u, x)v = e^{xL(-1)} Y(v, -x) u
$$
for $u, v \in V$, under the condition that $L(-1)I \subset I$ the
``left-ideal'' and ``right-ideal'' conditions (\ref{co1}) and
(\ref{co2}) are equivalent. In particular, (\ref{co1}) and (\ref{co2})
are equivalent for a vertex operator algebra (cf. Remark 3.9.8 in
\cite{LL}).
\end{remark}

\begin{proposition} \label{p}
The space $N^1(\Lambda_0)$ is the ideal of the vertex operator algebra
$N(\Lambda_0)$ generated by $x_{\alpha}(-1)^2 \cdot v_{\Lambda_0}^N$.
\end{proposition}

{\bf Proof:} By Remark 6.2.24 and Proposition 6.6.17 in \cite{LL} we
have that $N^1(\Lambda_0)$ is an ideal of $N(\Lambda_0)$. Since
$x_{\alpha}(-1)^2 \cdot v_{\Lambda_0}^N$ belongs to both
$N^1(\Lambda_0)$ and the ideal generated by $x_{\alpha}(-1)^2 \cdot
v_{\Lambda_0}^N$, and for any $a \in \goth{g}$ and $m \in \mathbb{Z}$,
$a(m)$ is a component of the vertex operator $Y(a(-1)\cdot
v_{\Lambda_0}^N, x)$, the statement follows. $\;\;\; \Box$
\vspace{1em}

As in ring theory, we call an ideal of a vertex (operator) algebra
generated by one element a {\it principal} ideal.  Thus Proposition
\ref{p} says that $N^1(\Lambda_0)$ is the principal ideal of
$N(\Lambda_0)$ generated by the ``null vector'' $x_{\alpha}(-1)^2
\cdot v_{\Lambda_0}^N$ (cf. (\ref{lie}) for the Lie-algebraic
statement).

The restrictions of the maps (\ref{Pi}) to the principal subspaces
$W^N(\Lambda_i)$ are the linear maps $\pi_{\Lambda_i}$ (\ref{pi})
introduced in Section 2. Thus
\begin{equation} \label{k}
\mbox{Ker} \; \pi_{\Lambda_i}= N^1(\Lambda_i) \cap W^N(\Lambda_i)
\end{equation} 
for $i=0,1$.

\begin{remark}
\rm We observe that $N^1(\Lambda_0) \cap W^N(\Lambda_0)$ is an ideal
of $W^N(\Lambda_0)$. Indeed, since $N^1(\Lambda_0)$ is an ideal and
$W^N(\Lambda_0)$ is a vertex subalgebra of $N(\Lambda_0)$,
$$
Y(v, x)w \in (N^1(\Lambda_0) \cap W^N(\Lambda_0)) ((x)) \; \; \mbox{and} \; \; 
Y(w, x)v \in (N^1(\Lambda_0) \cap W^N(\Lambda_0)) ((x))
$$
for $v \in W^N(\Lambda_0)$ and $w \in N^1(\Lambda_0) \cap
W^N(\Lambda_0)$.  Also, $N^1(\Lambda_1) \cap W^N(\Lambda_1)$ is a
$W^N(\Lambda_0)$-submodule of $W^N(\Lambda_1)$.
\end{remark} 

The intersection $N^1(\Lambda_0) \cap W^N(\Lambda_0)$, which equals
$I_{\Lambda_0} \cdot v_{\Lambda_0}^N$ by (\ref{k}) and Theorem
\ref{th2}, is also a principal ideal, generated by the same null vector:

\begin{proposition} \label{4}
We have that $I_{\Lambda_0} \cdot v_{\Lambda_0}^N$ is the ideal of the
vertex algebra $W^N(\Lambda_0)$ generated by $x_{\alpha}(-1)^2 \cdot
v_{\Lambda_0}^N$.  
\end{proposition}   

{\bf Proof:} We denote by $I$ the ideal of $W^N(\Lambda_0)$ generated
by $x_{\alpha}(-1)^2 \cdot v_{\Lambda_0}^N$.  By (\ref{ideal1}) we
have
\begin{equation} \label{new}
I_{\Lambda_0} \cdot v_{\Lambda_0}^N= \sum_{t \geq 2}
U(\bar{\goth{n}}_{-}) R_t^0 \cdot v_{\Lambda_0}^N,
\end{equation}
and this space contains $x_{\alpha}(-1)^2 \cdot v_{\Lambda_0}^N$.

We first show the inclusion 
\begin{equation} \label{egal1}
I_{\Lambda_0} \cdot v_{\Lambda_0}^N \subset I.
\end{equation}
Since $x_{\alpha}(-1)^2 \cdot v_{\Lambda_0}^N \in I$ and $I$ is an
ideal of $W^N(\Lambda_0)$ we have
$$Y(x_{\alpha}(-1)^2 \cdot v_{\Lambda_0}^N, x) v_{\Lambda_0}^N \in I((x)).$$
But since
\begin{equation} \label{quadratics}
Y(x_{\alpha}(-1)^2 \cdot v_{\Lambda_0}^N, x) = \sum_{t \in \mathbb{Z}}
\left ( \sum_{m_1+m_2=-t} x_{\alpha}(m_1)x_{\alpha}(m_2) \right )
x^{-t-2}
\end{equation}
we see that
\begin{equation} \label{s}
R_t^0 \cdot v_{\Lambda_0}^N \in I
\end{equation}
for each $t \geq 2$.  Finally, since each $x_{\alpha}(m)$, $m \in
\mathbb{Z}$, is a component of the vertex operator $Y(x_{\alpha}(-1)
\cdot v_{\Lambda_0}^N, x)$, the inclusion (\ref{egal1}) holds.

It remains only to show that $I_{\Lambda_0} \cdot v_{\Lambda_0}^N$ is
an ideal of $W^N(\Lambda_0)$.  We first observe that the operator
$L(-1)$ preserves $I_{\Lambda_0} \cdot v_{\Lambda_0}^N$.  Indeed, by
using (\ref{L(-1)}), (\ref{bracket}) and (\ref{quadratics}) we obtain
\begin{equation} \label{l(-1)r}
L(-1)(R_t^0 \cdot v_{\Lambda_0}^N)= (t-1) R_{t+1}^0 \cdot
v_{\Lambda_0}^N \in I_{\Lambda_0} \cdot v_{\Lambda_0}^N
\end{equation}
for any $t \geq 2$. More generally, from (\ref{b}), for any $m_1,
\dots, m_r \geq 1$ we have
\begin{eqnarray} \nonumber
&& L(-1)(x_{\alpha}(-m_1) \cdots x_{\alpha}(-m_r) R_t^0 \cdot
v_{\Lambda_0}^N) \nonumber \\
&& \hspace{1em} = \sum_{j=1}^r m_jx_{\alpha}(-m_1) \cdots
x_{\alpha}(-m_j-1) \cdots x_{\alpha}(-m_r) R_t^0 \cdot v_{\Lambda_0}^N
\nonumber \\
&& \hspace{3em} + (t-1)x_{\alpha}(-m_1) \cdots x_{\alpha}(-m_r)R_{t+1}^0
\cdot v_{\Lambda_0}^N \in I_{\Lambda_0} \cdot v_{\Lambda_0}^N,
\nonumber
\end{eqnarray}
and this shows that
\begin{equation} \label{l}
L(-1) (I_{\Lambda_0} \cdot v_{\Lambda_0}^N) \subset I_{\Lambda_0}
\cdot v_{\Lambda_0}^N
\end{equation}
(recall (\ref{new})).

In view of Remark \ref{skew} and (\ref{l}), in order to prove that
$I_{\Lambda_0} \cdot v_{\Lambda_0}^N$ is an ideal of $W^N(\Lambda_0)$
it is enough to show that $v_n w \in I_{\Lambda_0} \cdot
v_{\Lambda_0}^N$ for any $v \in W^N(\Lambda_0)=U(\bar{\goth{n}}_{-})
\cdot v_{\Lambda_0}^N$, $w \in I_{\Lambda_0} \cdot v_{\Lambda_0}^N$
and $n \in \mathbb{Z}$.  This will follow once we show that $v_n w \in
U(\bar{\goth{n}})w$.

To do this, note that $v$ is a linear combination of monomials in
$x_{\alpha}(-m)$, with each $m \geq 1$, applied to $v_{\Lambda_0}^N$.
We first assume that $v = v_{\Lambda_0}^N$. Then
\begin{equation}
Y(v,x)w= Y(v_{\Lambda_0}^N, x)w= w 
\end{equation}
(since $Y(v_{\Lambda_0}^N, x)$ is the identity operator), so that $v_n
w \in U(\bar{\goth{n}})w$.  Now we assume that $v$ is a linear
combination of elements $x_{\alpha}(-m_1) \cdots x_{\alpha}(-m_r)
\cdot v_{\Lambda_0}^N$ with $m_1, \dots, m_r \geq 1$, which we write
as $x_{\alpha}(-m_1) u$ for $u= x_{\alpha}(-m_2) \cdots
x_{\alpha}(-m_r) \cdot v_{\Lambda_0}^N$.  Since
$$
x_{\alpha}(-m_1)u= (x_{\alpha}(-1) \cdot v_{\Lambda_0}^N)_{-m_1}u,
$$
the iterate formula (cf. formula (3.1.11), and more specifically,
formula (3.8.12), in \cite{LL}) gives
\begin{eqnarray} \label{induc}
&& Y(v, x)w= Y \left ( (x_{\alpha}(-1) \cdot v_{\Lambda_0}^N)_{-m_1}
u, x \right ) w \\ 
&&= \frac{1}{(m_1-1)!} \left ( \left (\frac{d}{dx}
\right )^{m_1-1} Y(x_{\alpha}(-1) \cdot v_{\Lambda_0}^N, x) \right )
Y(u, x) w \nonumber \\ 
&&= \frac{1}{(m_1-1)!}  \left ( \left (
\frac{d}{dx} \right )^{m_1-1} \sum_{j \in {\mathbb{Z}}}x_{\alpha}(j)x^{-j-1}
\right ) Y(u, x) w.
\nonumber
\end{eqnarray}
Now (\ref{induc}) and induction on $r \geq 0$ imply that $v_n w \in
U(\bar{\goth{n}}) w$, completing the proof.  $\;\;\; \Box$
\vspace{1em}

Using the notation $(v)_V$ for the ideal generated by an element $v$
of a vertex (operator) algebra $V$, we have succeeded in reformulating
(\ref{eqone}) in Theorem \ref{th2} as follows:
\begin{theorem} \label{th3}
We have
\begin{equation}
{\rm{Ker}}\; \pi_{\Lambda_0}= (x_{\alpha}(-1)^2 \cdot v_{\Lambda_0}^N
) _{N(\Lambda_0)} \cap W^N(\Lambda_0)= (x_{\alpha}(-1)^2 \cdot
v_{\Lambda_0}^N)_{W^N(\Lambda_0)}.
\end{equation}
In particular, the intersection with the vertex subalgebra
$W^N(\Lambda_0)$ of the principal ideal of $N(\Lambda_0)$ generated by
the null vector $x_{\alpha}(-1) ^2 \cdot v_{\Lambda_0}^N$ coincides
with the principal ideal of the vertex subalgebra $W^N(\Lambda_0)$
generated by the same null vector.
\end{theorem}

\noindent {\small \sc Department of Mathematics, Rutgers University,
Piscataway, NJ 08854} \\
Current address:\\
\noindent {\small \sc Department of Mathematics, Ohio State University,
Columbus, OH 43210} \\ 
{\em E--mail address}:
calinescu@math.ohio-state.edu\\
On leave from the Institute of Mathematics of the Romanian Academy.\\
 
\vspace{2mm}
\noindent {\small \sc Department of Mathematics, Rutgers University,
Piscataway, NJ 08854} \\ 
{\em E--mail address}:
lepowsky@math.rutgers.edu \\

\vspace{2mm}
\noindent {\small \sc Department of Mathematics and Statistics,
University at Albany (SUNY), Albany, NY 12222} \\ 
{\em E--mail
address}: amilas@math.albany.edu

\end{document}